\documentclass[11pt, a4paper, reqno]{amsart}
\usepackage[a4paper, margin=3.75cm]{geometry}
\usepackage[T2A]{fontenc}
\usepackage[cp866]{inputenc}
\usepackage[russian,english]{babel}
\usepackage{amsfonts,amssymb,mathrsfs,amscd,comment}
\usepackage{amsthm,graphicx} 
\def\NoBlackBoxes{\overfullrule0pt}

\NoBlackBoxes
\theoremstyle{plain}
\newtheorem{theorem}{Theorem}

\theoremstyle{definition}

\usepackage[breaklinks]{hyperref}
\hypersetup{
unicode=true,
colorlinks=true,
linkcolor=blue,
citecolor=blue,
urlcolor=blue,
filecolor=blue,
bookmarksnumbered=true,
pdfstartview=FitH,
pdfhighlight=/N
}
\theoremstyle{main}
\let\savedef=\endproof
\def\endproof{~$\square$\savedef}
\def\bad{\spaceskip=0.33emplus0.6emminus0.15em\immediate\write5{\string\bad}}

\theoremstyle{plain}

\let\myh\widehat

\let\eps\varepsilon

\def\HH{\mathscr H}
\def\MM{\mathscr M}

\def\CC{\mathbb C}

\def\NN{\mathbb N}

\def\mdeg{\operatorname{deg}}

\def\mcap{\operatorname{cap}}

\def\supp{\operatorname{supp}}
\def\const{\operatorname{const}}
\def\<{\left\langle}
\def\>{\right\rangle}
\def\({\left(}
\def\){\right)}
\def\[{\left[}
\def\]{\right]}
\NoBlackBoxes

\def\bad{\spaceskip=0.33emplus0.6emminus0.15em\immediate\write5{\string\bad}}
\def\NN{\mathbb N}

\def\CC{\mathbb C}
\def\PP{\mathbb P}

\def\zz{\mathbf z}
\def\RS{\mathfrak R}

\def\dzeta{\mathbf\zeta}

\def\bGamma{\boldsymbol\Gamma}

\let\myh\widehat
\let\myt\widetilde

\def\const{\operatorname{const}}

\def\({\left(}
\def\){\right)}
\def\[{\left[}
\def\]{\right]}
\def\<{\left\langle}
\def\>{\right\rangle}

\let\geq\geqslant

\let\leq\leqslant
\begin{document}

\selectlanguage{english}

\title{A direct proof of Stahl's theorem for a generic class of algebraic functions}
\author[Sergey~P.~Suetin]{Sergey~P.~Suetin}
\address{Steklov Mathematical Institute of the Russian Academy of Sciences, Russia}
\email{suetin@mi-ras.ru}

\date{July 1, 2021}

\maketitle

\markright{On a Stahl's theorem}

\begin{abstract}
Under the assumption of the existence of Stahl's $S$-compact set we give a short proof of the limit zeros distribution of Pad\'e polynomials and convergence in capacity of diagonal Pad\'e approximants for a generic class of algebraic functions. The proof is direct but not from the opposite as Stahl's original proof is. The generic class means in particular that all branch points of the multi-sheeted Riemann surface of the algebraic function are of the first order (i.e., we assume the surface is such that all branch points are of square root type).

We do not use the relations of orthogonality at all. The proof is based on the maximum principle only.

Bibliography:~\cite{Sta12}~titles.
\end{abstract}

\section{Introduction}\label{s1}

\subsection{}\label{s1s1}
The seminal Stahl's Theory (see~\cite{Sta97b},~\cite{Sta12},~\cite{ApBuMaSu11} and the bibliography therein) on the convergence of diagonal Pad\'e approximants for the multivalued analytic functions consists of two parts. Namely, geometrical part and analytical part. In the first part it is proven that for a multi-valued function there exists a unique admissible compact set $S$ possesses Stahl's symmetric property. Traditionally it is refereed to as an $S$-compact set or briefly $S$-curve (see~\cite{Sta97b} and~\cite{Rak12} for details). In the original Stahl's Theory it is supposed that the set of singularities of a multivalued analytic function is of zero (logarithmic) capacity. Here we restrict our investigation to the case of finite number of singular points only. Note that a short proof of the existence of an $S$-curve for a multivalued analytic function with a finite number of singular points was proposed by E.~Rakhmanov in 1994 in the unpublished paper~\cite{PeRa94} (see also~\cite{MaRaSu11}). Rakhmanov's proof is based on the connection between capacity and potential of equilibrium measure of a compact set.

In the second part of the Stahl's Theory it is proven that there exists zeros limit distribution of Pad\'e polynomials and convergence in capacity of the diagonal Pad\'e approximants takes place as well. The original Stahl's proof of the existence of limit zeros distribution of Pad\'e polynomials is based on the method of proof from the opposite. In the current  paper we give a short and direct proof of the second part of Stahl's Theory for a generic class of algebraic functions and under the assumption that the first geometrical component is valid.
We do not use the relations of orthogonality at all.
The proof is based on the maximum principle only (cf.~\cite{Nut84} and~\cite{Sta87}).

Let introduce a generic class $\mathscr F$ of admissible multivalued functions which was mentioned above.

We set that $f\in\mathscr F$ if the following conditions are satisfied:

(1) $f$ is an algebraic function with all the branch points of
the first order (i.e., we assume that all the branch points are of square root type);

(2) $S$-curve of $f$ consists of a finite number of non-intersected closed analytic arcs each of which contains exactly two branch points of $f$ -- the ends of the arc.

Some discussion on this assumptions are given in Section~\ref{s3}.

\subsection{}\label{s1s2}
Let $f\in\mathscr F$, $f_\infty\in\HH(\infty)$ and $S=S(f_\infty)$ be the corresponding Stahl's compact set.
Let $\lambda_S$ be the probability equilibrium measure supported on $S$ and $V^{\lambda_S}(z)$ be the corresponding logarithmic potential of $\lambda_S$:
\begin{equation}
V^{\lambda_S}(z)=\gamma_S-g^{}_S(z,\infty), \quad z\in D:=\myh{\mathbb C}\setminus S,
\label{3}
\end{equation}
$g^{}_S(z,\infty)$ is the Green function for the domain $D$ with the logarithmic singularity at the infinity point $z=\infty$ and $\gamma_S$ is the Robin constant for $D$.

Set $\PP_n:=\mathbb C_n[z]$ be the set of all algebraic polynomials with complex coefficients of order $\leq{n}$. For arbitrary polynomial $Q\in\CC[z]$, $Q\not\equiv0$, we denote by $\chi(Q)$ the zero counting measure of $Q$: $\chi(Q):=\sum\limits_{\zeta:Q(\zeta)=0}\delta_\zeta$, where each zero of $Q$ is represented with respect to its multiplicity.

For $f_\infty\in\HH(\infty)$ and arbitrary $n\in\NN$ \ Pad\'e polynomials $P_n,Q_n\in\PP_n$, $Q_n\not\equiv0$,  are defined (non-uniquely) from the following relation:
\begin{equation}
R_n(z):=(Q_nf_\infty-P_n)(z)=O\(\frac1{z^{n+1}}\),\quad z\to\infty;
\label{4}
\end{equation}
$R_n(z)$ is the so-called error function.

The following result holds (see~\cite{Sta97b},~\cite{Sta12}).

\begin{theorem}\label{the1}
Let $f\in\mathscr F$ and $f_\infty\in\HH(\infty)$. Then as $n\to\infty$
\begin{equation}
\frac1n\chi(Q_n),\frac1n\chi(P_n)\overset{*}\longrightarrow\lambda_S.
\label{5}
\end{equation}
Set $Q_n(z)=z^{k_n}+\dotsb$, then as $n\to\infty$
\begin{gather}
|R_n(z)|^{1/n}\overset{\mcap}\longrightarrow \mcap(S)e^{-g^{}_S(z,\infty)}\quad\text{inside $D$},
\label{6}\\
\bigl|f(z)-[n/n]_f(z)\bigr|^{1/n}
\overset{\mcap}\longrightarrow e^{-2g^{}_S(z,\infty)}\quad\text{inside $D$}.
\label{7}
\end{gather}
\end{theorem}

Recall that the existence of an $S$-compact set for a function $f\in\mathscr F$ follows from~\cite{PeRa94}.

\section{Proof of Theorem~\ref{the1}}\label{s2}

\subsection{}\label{s2s1}
For two positive sequences $\{\alpha_n\}$ and $\{\beta_n\}$ the relation $\alpha_n\asymp \beta_n$ means that  $0<C_1\leq \alpha_n/\beta_n\leq C_2<\infty$ for $n=1,2,\dots$ and some constants $C_1,C_2$ which do not depend on $n$. For two sequences $\{\alpha_n(z)\}$ and $\{\beta_n(z)\}$ of functions holomorphic in a domain $\Omega$ the relation $\alpha_n\asymp\beta_n$ means that for each compact set $K\subset\Omega$ and $n=1,2,\dots$ \  $0<C_1\leq |\alpha_n(z)/\beta_n(z)|\leq C_2<\infty$ for $z\in K$ where the constants $C_1,C_2$ depend on $K$  but do not depend on $n$ and $z\in K$. Evidently that for such pairs of sequences and functions we have that $|\alpha_n/\beta_n|^{1/n}\to1$ as $n\to\infty$.

Since $f\in\mathscr F$, we have that $S=S(f)=\bigsqcup\limits_{j=1}^p S_j$, where $S_j=\operatorname{arc}(a_{2j-1},a_{2j})$.
Set $w^2=\prod\limits_{j=1}^p(z-a_{2j-1})(z-a_{2j})$. Then the two-sheeted hyperelliptic Riemann surface (RS) $\RS_2(w)$ of the function $w$ is the RS
 associated with $f_\infty\in\HH(\infty)$, $f\in\mathscr F$, in accordance with Stahl's Theory. Point $\zz$ on $\RS_2(w)$ is given by $\zz=(z,w)$. The RS $\RS_2(w)$ may be considered as a two-sheeted covering of the Riemann sphere $\myh{\mathbb C}$ with the canonical projection $\pi$, $\pi\colon\RS_2(w)\to\myh{\mathbb C}$, $\pi(\zz)=z$.  Let $\bGamma:=\pi^{-1}(S)$. Then $\Gamma$ shares $\RS_2(w)$ into two domains, We shall call them the (open) sheets of $\RS_2(w)$.
 The function $w$ is a single-valued function on this RS and has opposite signs on the sheets. We denote by $\RS_2^{(0)}$ the sheet of $\RS_2(w)$ where $w(z)/z^{p}\to1$ as $z\to\infty$ and refereed to it as zero sheet. Another sheet will be denoted by $\RS^{(1)}_2$ and refereed to as first sheet of $\RS_2(w)$. Then we have $\RS_2(w)=\RS_2^{(0)}\sqcup\Gamma\sqcup\RS_2^{(1)}$. Points on the sheets will be denoted by $z^{(j)}$, $j=0,1$, respectively. Evidently $\pi(\RS_2^{(j)})=D$. We shall identify $\RS_2^{(0)}$ with the Stahl's domain $D=\myh{\mathbb C}\setminus{S}$ and $\infty^{(0)}$ with $\infty$ and consider the germ $f_\infty$ as the germ $f_{\infty^{(0)}}$. In general $\RS_2(w)$ is not the Riemann surface of $f$ and thus $f$ is not a single-valued on $\RS_2(w)$.  But since $f\in\mathscr F$ the germ $f_{\infty^{(0)}}$ can be extended from the infinity point $\infty^{(0)}$ to the whole zero sheet $\RS^{(0)}_2$  and even more to a small enough neighborhood $V^{(0,1)}$ of  $\bGamma$ , $V^{(0,1)}\cap\RS_2^{(1)}\neq\varnothing$. In accordance with Stahl's Theory (see also~\cite{PeRa94}) the Green function $g^{}_S(z,\infty)$ of $D$  can be lifted to RS $\RS_2(w)$ as $g(\zz)$ with the following properties: $g(z^{(0)})=g(z,\infty)$, $g(z^{(1)})=-g(z^{(0})<0$.  From now on we shall suppose that $V^{(0,1)}$ is such that $\partial V^{(0,1)}\cap\RS_2^{(1)}=\{z^{(1)}:g(z,\infty)=\log{R}\}$, $R>1$. Set
$\mathfrak D:=\RS_2^{(0)}\cup V^{(0,1)}$ be a domain on $\RS_2(w)$ and $f_{\infty^{(0)}}$ continues in $\mathfrak D$ as a meromorphic  (analytic and single-valued) function, $f\in\MM(\mathfrak D)$.
The function $R_n(z)$ is also extended to $\mathfrak D$ as a meromorphic function $R_n(\zz)$.
Let a polynomial $q_m(z)=z^m+\dotsb$ is such that the set of its zeros coincide with the projection of poles of $f$ in $\mathfrak D$. Then $\myt{f}:=q_nf$ and $q_mR_n$ are holomorphic functions in $\mathfrak D$.

\subsection{}\label{s2s2}
For each $\rho\in(1,R)$ we denote by $\Gamma^{(1)}_\rho$ the set of points $z^{(1)}$ such that $g^{}_S(z,\infty)=\log\rho$ as $z^{(1)}\in\Gamma^{(1)}_\rho$. Clearly we have $g(z^{(1)})=-\log\rho$.  The set $\Gamma^{(0)}_\rho$ is defined similarly and we have on that set $g(z^{(0)})=\log\rho$.
For $\rho\in(1,R)$ let $D^{(1)}_\rho$ be a subdomain of $\mathfrak D$ with $\partial D^{(1)}_\rho=\Gamma^{(1)}_\rho$, $\infty^{(0)}\in D^{(1)}$. Similarly $D^{(0)}_\rho\subset\mathfrak D$, $\partial D^{(0)}=\Gamma^{(0)}_\rho$, $\infty^{(0)}\in D^{(0)}$. Set
\begin{equation}
u_n(\zz):=\log|q_m(z)R_n(\zz)|+(n+1-m)g(\zz), \quad \zz\in\mathfrak D.
\label{8}
\end{equation}
Function $u_n$ is a subharmonic function in $\mathfrak D$. Thus for each $\rho\in(1,R)$ we have that
\begin{equation}
u_n(\zz)\leq\max_{\dzeta\in\Gamma^{(1)}_\rho}u_n(\zeta),
\quad \zz\in D^{(1)}_\rho.
\label{9}
\end{equation}
From~\eqref{9} it follows that
\begin{equation}
\bigl|R_n(z^{(0)})q_m(z)\bigr|\leq\frac1{\rho^{2n+2-2m}}
M_{n,1}(\rho),
\qquad z\in\Gamma_\rho.
\label{10}
\end{equation}
where we set $M_{n,1}(\rho):=\max\limits_{z\in\Gamma_\rho}|R_n(z^{(1)})q_m(z)|$.
Relations~\eqref{8} and~\eqref{9} imply that for $1<\rho_1<\rho_2<R$ we have
\begin{equation}
M_{n,1}(\rho_1)\leq\(\frac{\rho_2}{\rho_1}\)^{n-m+1}M_{n,1}(\rho_2).
\label{11}
\end{equation}

The following identity holds true:
\begin{align}
R_n(z^{(0)})&=Q_n(z)f(z^{(0)})-P_n(z)
=Q_n(z)f(z^{(1)})-P_n(z)+Q_n(z)[f(z^{(0)})-f(z^{(1)})]\notag\\
&=R_n(z^{(1)})+Q_n(z)[f(z^{(0)})-f(z^{(1)})],\quad z\in \Gamma_\rho.
\label{12}
\end{align}
From~\eqref{10} and~\eqref{11} it follows that
\begin{equation}
\max_{z\in \Gamma_\rho}|Q_n(z)[f(z^{(0)})-f(z^{(1)})]q_m(z)|
=M_{n,1}(\rho)(1+\eps_n),\quad \eps_n\to0.
\label{13}
\end{equation}
From now on we shall consider only such $\rho\in(1,R)$ that $|q_m(z)[f(z^{(0)})-f(z^{(1)})]|\geq C(\rho)>0$ as $z\in\Gamma_\rho$ (clearly, $f(z^{(0)})-f(z^{(1)})\not\equiv0$).
Set $m_n(\rho):=\max\limits_{z\in \Gamma_\rho}|Q_n(z)|$. Then from~\eqref{13} it follows that
\begin{equation}
m_n(\rho)\asymp M_{n,1}(\rho).
\label{14}
\end{equation}

Since $\mdeg{Q_n}\leq{n}$, then Bernstein--Walsh theorem gives us the inequality:
\begin{equation}
|Q_n(z)|\leq e^{n g^{}_{\Gamma_{\rho_1}}(z,\infty)}m_n(\rho_1),\quad z\in\Gamma_{\rho_2},\quad \rho_2>\rho_1,
\label{15}
\end{equation}
where $g^{}_{\Gamma_{\rho_1}}(z,\infty)$ is the Green function for the domain $g^{}_S(z,\infty)>\log\rho_1$. Clearly, $g^{}_{\Gamma_{\rho_1}}(z,\infty)=g^{}_S(z,\infty)-\log\rho_1$. From~\eqref{15} we obtain the estimate
\begin{equation}
m_n(\rho_2)\leq \(\frac{\rho_2}{\rho_1}\)^{n}m_n(\rho_1).
\label{16}
\end{equation}
After combining the relations~\eqref{11},~\eqref{14} and~\eqref{16} we obtain that
\begin{equation}
m_n(\rho_2)\asymp \(\frac{\rho_2}{\rho_1}\)^n m_n(\rho_1)\quad
\text{and} \quad
M_{n,1}(\rho_2)\asymp \(\frac{\rho_2}{\rho_1}\)^n M_{n,1}(\rho_1).
\label{17}
\end{equation}

Let $Q_n(z)=z^{k_n}+\dotsb$, $k_n=\mdeg{Q_n}\leq{n}$.
Again Bernstein--Walsh theorem gives us the inequality:
\begin{equation}
|Q_n(z)|\leq e^{k_n g^{}_{\Gamma_{\rho_1}}(z,\infty)}m_n(\rho_1),
\quad z\in\Gamma_{\rho_2},\quad \rho_2>\rho_1.
\label{18}
\end{equation}
From this relation it follows that
\begin{equation}
m_n(\rho_2)\leq\(\frac{\rho_2}{\rho_1}\)^{k_n}m_n(\rho_1).
\label{19}
\end{equation}
But in accordance with~\eqref{17} $m_n(\rho_2)\asymp (\rho_2/\rho_1)^n m_n(\rho_1)$. From here and~\eqref{19} we obtain the inequality
\begin{equation}
\(\frac{\rho_2}{\rho_1}\)^{n}\leq C\(\frac{\rho_2}{\rho_1}\)^{k_n},
\label{20}
\end{equation}
where $C=C(\rho_1,\rho_2)$, $k_n\leq{n}$, $1<\rho_1<\rho_2$. It directly follows from~\eqref{20} that
\begin{equation}
\mdeg{Q_n}/n\to1,\quad n\to\infty.
\label{21}
\end{equation}

\subsection{}\label{s2s3}
Let $D_\rho:=\{z\in\myh{\mathbb C}: g^{}_S(z,)>\log\rho\}$, $\rho>1$,  $g^{}_{\Gamma_\rho}(z,\infty)$  be the Green function for $D_\rho$. Then $g^{}_{\Gamma_\rho}(z,\infty)=g^{}_S(z,\infty)-\log\rho=\log{|z|}+\gamma_\rho+o(1)=\log|z|+\gamma_S-\log\rho+o(1)$. Thus $\gamma_\rho=\gamma_S-\log\rho$.

Set
\begin{equation}
u_n(z):=\frac1{k_n}\log|Q_n(z)|-g^{}_{\Gamma_\rho}(z,\infty)-\frac1{k_n}\log m_n(\rho).
\label{22}
\end{equation}
Since function $u_n$ is subharmonic in the domain $D_\rho$  and $u_n\leq0$ on $\Gamma_\rho$, then $u_n\leq0$ in $D_\rho$ and $u_n(\infty)\leq0$. Thus we have that $\log\rho-\gamma_S\leq\log m_n(\rho)^{1/k_n}$ and finally obtain
\begin{equation}
m_n(\rho)^{1/k_n}\geq  \rho e^{-\gamma_S}=\rho\mcap(S).
\label{23}
\end{equation}

For $z\in K\Subset D\rho$ and $\zeta\in\Gamma_\rho$ we have by definition $g^{}_{\Gamma_\rho}(z,\zeta)=0$. Let now extend the function $g^{}_{\Gamma_\rho}(z,\zeta)$ inside $\Gamma_\rho$ by identity $g^{}_{\Gamma_\rho}(z,\zeta)\equiv0$ when $z\in D_\rho$, $\zeta \in\operatorname{int} \Gamma_\rho$.

Let $Q_n(z)=\prod\limits_{j=1}^{k_n}(z-\zeta_j)$,
\begin{equation}
v_n(z):=\frac1{k_n}\log|Q_n(z)|-g^{}_{\Gamma_\rho}(z,\infty)+\frac1{k_n}\sum_{j=1}^{k_n} g^{}_{\Gamma_\rho}(z,\zeta_j)-\frac1{k_n}\log m_n(\rho),
\label{24}
\end{equation}
$z\in D_{\rho}$. Then $v_n$ is a harmonic function in $D_\rho$ and $v_n\leq0$ on $\Gamma_\rho$. From this it follows that for $z\in \Gamma_{\rho_2}$, $\rho<\rho_2<R$,  we have
\begin{equation}
|Q_n(z)|\exp\left\{\sum_{j=1}^{k_n} g^{}_{\Gamma_\rho}(z,\zeta_j)\right\}
\leq m_n(\rho)
 \exp\biggl\{k_ng^{}_{\Gamma_\rho}(z,\infty) \biggr\}.
\label{25}
\end{equation}
Let $z_n^*\in \Gamma_{\rho_2}$ is such that $|Q_n(z_n^*)|=m_n(\rho_2)$. Then from~\eqref{25} we obtain that
\begin{equation}
m_n(\rho_2)\exp\biggl\{\sum_{j=1}^{k_n}g^{}_{\Gamma_\rho}(z_n^*,\zeta_j)  \biggr\}\leq m_n(\rho)\(\frac{\rho_2}{\rho}\)^{k_n}.
\label{26}
\end{equation}
Ultimately from~\eqref{26} and~\eqref{17} we obtain that
\begin{equation}
\exp\biggl\{\frac1{k_n}\sum_{j=1}^{k_n}g^{}_{\Gamma_\rho}(z_n^*,\zeta_j) \biggr\}\leq C(\rho,\rho_2)^{1/k_n},\quad \eps_n\to0.
\label{27}
\end{equation}
Let now $\Lambda\subset\NN$ is a subsequence such that $z_n^*\to z^*\in\Gamma_{\rho_2}$ and $\frac1{k_n}\chi(Q_n)\to \mu$ as $n\to\infty$, $n\in\Lambda$, $\mu(1)=1$.
Then according to the descendence principle (see~\cite[Chapter I, \S~3, Theorem 1.3]{Lan66},~\cite{SaTo97},~\cite{Chi18} and~\cite{Chi20}), we have that
\begin{equation}
\int g^{}_{\Gamma_\rho}(z^*,\zeta)\,d\mu(\zeta)\leq 0,
\quad z^*\in\Gamma_{\rho_2},\quad \rho_2>\rho.
\label{28}
\end{equation}
It directly follows from~\eqref{28} that $\supp{\mu}\subset\myh{\mathbb C}\setminus D_\rho$. Since it is true for every $\rho>1$ than we have that $\supp{\mu}\subset S$.

\subsection{}\label{s2s4}
So we have got that $k_n/n\to1$ and each limit point $\mu$ of the sequence $\{\frac1{k_n}\chi(Q_n)\}$ satisfies the condition $\supp\mu\subset S$.

Since $\supp{\mu}\subset{S}$ we have that
$$
\varlimsup_{\substack{n\to\infty\\n\in\Lambda}}
m_n(\rho)^{1/k_n}\leq C_2<\infty.
$$
Also from~\eqref{23} it follows that
$$
\varliminf_{\substack{n\to\infty\\n\in\Lambda}}
m_n(\rho)^{1/k_n}\geq C_1>0.
$$
Set
\begin{equation}
u_n(z):=\frac1{k_n}\log\frac1{|Q_n(z)|}-\frac1{k_n}\sum_{j=1}^{k_n}
g^{}_{\Gamma_\rho}(z,\zeta_j)+g^{}_{\Gamma_\rho}(z,\infty),
\quad z\in D_\rho.
\notag
\end{equation}
Then $\{u_n\}$ is a sequence of harmonic functions in $D_\rho$.
Let $u(z)=\lim\limits_{n\to\infty}u_n(z)$. Then $u$ is a harmonic function with the following properties (see~\eqref{17}):
\begin{gather}
u(z)=V^{\mu}(z)-V^{\lambda_S}(z)+\const,\notag\\
\min_{z\in\Gamma_{\rho_1}}u(z)=\min_{z\in\Gamma_{\rho_2}}u(z)
\quad\text{for each}\quad \rho_1,\rho_2>\rho.
\notag
\end{gather}
From this it follows that $u(z)=\const$ and thus $V^\mu(z)=V^{\lambda_S}(z)$ for $z\in D_\rho$. Since $\supp\mu\subset S$, then $V^\mu(z)=V^{\lambda_S}$ for $z\in D$. Finally $\mu=\lambda_S$ because $S$ has no inner points.

Eventually we have that for each $\rho>1$ there exist the limits
\begin{align}
\lim_{n\to\infty}m_n(\rho)^{1/n}=\lim_{n\to\infty}M_{n,1}(\rho)^{1/n}=\rho
\mcap{S},
\quad \rho>1,
\label{30}\\
\lim_{n\to\infty}\max_{z\in S}|Q_n(z)|^{1/n}=\mcap{S}=e^{-\gamma_S},\quad
\lim_{n\to\infty}\max_{z\in\Gamma_\rho}|R_n(z^{(0)})|^{1/n}=\frac1\rho
e^{-\gamma_S}.
\label{31}
\end{align}
From the above it easy to prove that
\begin{equation}
|R_n(z)|^{1/n}\overset\mcap\longrightarrow \mcap(S) e^{-g^{}_S(z,\infty)},
\quad n\to\infty.
\label{32}
\end{equation}
Indeed, let
\begin{equation}
v_n(z):=\frac1{n-m+1}\log|q_m(z)R_n(z)|+g^{}_S(z,\infty)+\gamma_S.
\label{33}
\end{equation}
Then for each $\rho>1$ function $v_n$ is a subharmonic function in $D_\rho$ and  from~\eqref{31} it follows that for $z\in\Gamma_\rho$ we have
\begin{equation}
v_n(z)\leq C_n,
\label{34}
\end{equation}
where $C_n\to0$ as $n\to\infty$. From~\eqref{34} we obtain that for $z\in D_\rho$
\begin{equation}
|q_m(z)R_m(z)|^{1/(n-m+1)}\leq\mcap(S)e^{-g^{}_S(z,\infty)}e^{C_n}.
\label{35}
\end{equation}
Thus
\begin{equation}
\varlimsup_{n\to\infty}|R_n(z)|^{1/n}\leq\mcap(S)e^{-g^{}_S(z,\infty)}.
\label{36}
\end{equation}
Relation~\eqref{32} follows from~\eqref{31} and~\eqref{36} by two-constants theorem (cf.~\cite[\S~3, subsection 8, (31)--(36)]{GoRa86}).

Since $\frac1n\chi(Q_n)\to\lambda_S$, then inside $D$
\begin{equation}
|Q_n(z)|^{1/n}\overset\mcap\longrightarrow e^{-V^{\lambda_S}(z)},
\quad n\to\infty.
\label{37}
\end{equation}
From~\eqref{32} and~\eqref{37} the relation~\eqref{7} follows.

Theorem~\ref{the1} is proven.

\section{Appendix}\label{s3}

In fact the admissible class of analytic functions is wider then it is described in the items (1)--(2) of definition of $\mathscr F$.
In particular our approach is valid for the class of multivalued functions  generated by the inverse Zhoukovskii transform. More precisely, let $\varphi(z):=z+(z^2-1)^{1/2}$, where $z\in\myh{\mathbb C}\setminus\Delta$, $\Delta=[-1,1]$ and such  branch of the function $(\cdot)^{1/2}$ is chosen that $\varphi(z)/z\to2$ as $z\to\infty$. Let $1<A<B<\infty$ and $a:=(A+1/A)/2$, $b:=(B+1/B)/2$.  Then the function
\begin{equation}
\mathfrak f_\Delta(z):=\[\(A-\frac1{\varphi(z)}\)\(B-\frac1{\varphi(z)}\)\]^{1/2}
\label{1}
\end{equation}
is an algebraic function of forth order with the set of square root singularities $\Sigma(\mathfrak f)=\{\pm1,a,b\}$. Under the above conditions on $(\cdot)^{1/2}$ we have that $\mathfrak f\in\HH(\myh{\mathbb C}\setminus\Delta)$ and Stahl's compact set $S(\mathfrak f)=[-1,1]$. Now let $\varphi_{\Delta_j}(z)$ be the inverse Zhoukovskii function for a segment $\Delta_j:=[\alpha_j,\beta_j]$, $j=1,\dots,m$, $\Delta_j\cap \Delta_k=\varnothing$, $j\neq k$.
Set
\begin{equation}
\mathfrak f(z):=\prod_{j=1}^m {\mathfrak f}_{\Delta_j}(z),
\label{2}
\end{equation}
where each $\mathfrak f_{\Delta_j}(z)$ is defined by~\eqref{1} with some suitable $A_j$ and $B_j$.

When all $\Delta_j\subset\mathbb R$, $j=1,\dots,m$, then $S(\mathfrak f)=\bigsqcup\limits_{j=1}^m\Delta_j$, Since $\mathbb C(z,\mathfrak f)\subset\mathscr F$, our approach is valid for each $f\in\mathbb C(z,\mathfrak f)$.

Another admissible and non-trivial subclass of $\mathscr F$ is obtained when at least one of the branch points $\alpha_j,\beta_j$, $j=1,\dots,m$, of $\mathfrak f$ does not belong to the real line.

We can generalize~\eqref{1} in the following way. Set
\begin{equation}
{\mathfrak f}_\Delta(z):=\prod\(A-\frac1{\varphi(z)} \)^\alpha\(B-\frac1{\varphi(z)}\)^\beta\cdot \dotsc\cdot\(C-\frac1{\varphi(z)}\)^\gamma,
\label{41}
\end{equation}
where $\alpha,\beta,\dots,\gamma\in\mathbb C\setminus\mathbb Z$, $\alpha+\beta+\dots+\gamma\in\mathbb Z$, to obtain functions which are not from $\mathscr F$ but for which the proof from Section 2 is also valid.


\def\by#1;{#1\unskip,}
\def\paper#1;{``#1\unskip''\unskip,}
\def\paperinfo#1;{#1\unskip.}
\def\eprint#1;{``#1\unskip''\unskip,}
\def\eprintinfo#1;{#1\unskip,}
\def\book#1;{``#1\unskip''\unskip,}
\def\inbook#1;{``#1\unskip''\unskip,}
\def\bookinfo#1;{#1\unskip,}
\def\jour#1;{#1\unskip,}
\def\issue#1;{#1\unskip,}
\def\yr#1;{#1\unskip,}
\def\pages#1.{#1\unskip.}
\def\vol#1;{\textbf{#1}\unskip,}
\def\finalinfo#1;{#1\unskip.}
\def\publ#1;{#1\unskip,}
\def\publadrr#1;{#1\unskip,}
\def\procinfo#1;{#1\unskip,}
\def\serial#1;{#1\unskip,}
\def\ed#1;{ed #1\unskip,}
\def\eds#1;{ed #1\unskip,}

\end{document}